\documentclass{article}
\usepackage{amsmath}
\usepackage{amssymb}
\usepackage{latexsym}
\usepackage{verbatim}
\usepackage{graphics}
\usepackage{amsfonts}
\usepackage{graphicx}
\usepackage{epstopdf}
\newtheorem{thm}{Theorem}
\newtheorem{conj}{Conjecture}
\newtheorem{prop}{Proposition}

\newtheorem{cor}{Corollary}
\newtheorem{quest}{Question}

\begin{document}

\title{Virtual Spatial Graphs}

\author{Thomas Fleming\\
   	Department of Mathematics\\
   	University of California, San Diego\\
   	La Jolla, CA 92093-0112\\
   	{\it  tfleming@math.ucsd.edu}
   \and
   	Blake Mellor\\
   	Mathematics Department\\
   	Loyola Marymount University\\
   	Los Angeles, CA  90045-2659\\
   	{\it  bmellor@lmu.edu}}

\date{}

\maketitle

\begin{abstract}

Two natural generalizations of knot theory are the study of spatially embedded graphs, and Kauffman's theory of virtual knots.  In this paper we combine these approaches to begin the study of virtual spatial graphs.

\end{abstract}

\tableofcontents

\section{Introduction} \label{S:intro}

Knot theory is the study of isotopy classes of circles (or, for links, disjoint unions of circles) embedded in 3-space.  There are many ways to extend the ideas of knot theory - two natural choices are the study of spatial graphs and Kauffman's theory of virtual knots \cite{ka}.  The theory of spatial graphs studies isotopy classes of general graphs embedded in 3-space; in particular, there has been considerable work done on spatial $\theta$-graphs \cite{go, nikk}.  Kauffman's theory of virtual knots goes in a very different direction.  Any knot can be described by its {\it diagram}, the result of projecting the embedding to a plane, retaining information about over- and under-crossings.  Such a projection can be described by its {\it Gauss code} - the sequence of crossings as we move around the knot.  However, there are many more such sequences than there are real knots; the problem of recognizing ``realizable" Gauss codes is an old one \cite{fo, de, rr}.  One motivation for virtual knots is to provide ``realizations" for the sequences which are not Gauss codes for classical 
knots. 

It is natural to combine these two generalizations.  In previous work, the authors extended the notion of Gauss codes to spatial graphs and looked at which codes were realizable by classical spatial graphs.  The goal of this paper is to begin the study of virtual spatial graphs, which provide a way to represent the ``non-realizable" Gauss codes.  We will give the basic definitions, a few fundamental properties, and provide some examples.  Future papers will continue various aspects of this study, such as looking at the notion of intrinsic linking in virtual spatial graphs.

\medskip
\noindent{\sc Acknowledgements:}  The authors would like to acknowledge the hospitality of Waseda University, Tokyo, and Professor Kouki Taniyama during the International Workshop on Knots and Links in a Spatial Graph in July, 2004, where the idea for this project was conceived.  During the research for this paper, the second author was supported by an LMU Faculty Research grant and Junior Faculty Sabbatical.  We would also like to thank V. Manturov and J. Uhing for valuable comments on an early draft of this paper.

\section{Defining Virtual Spatial Graphs} \label{S:definition}

\subsection{Definition and Reidemeister moves} \label{SS:reidemeister}

\begin{figure} [ht]
    $$\includegraphics{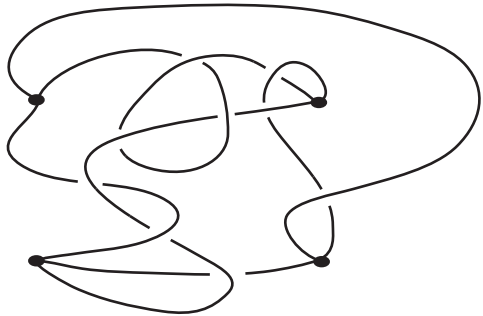}$$
    \caption{A graph diagram} \label{F:diagram}
    \end{figure}

\begin{figure} 
    $$\includegraphics{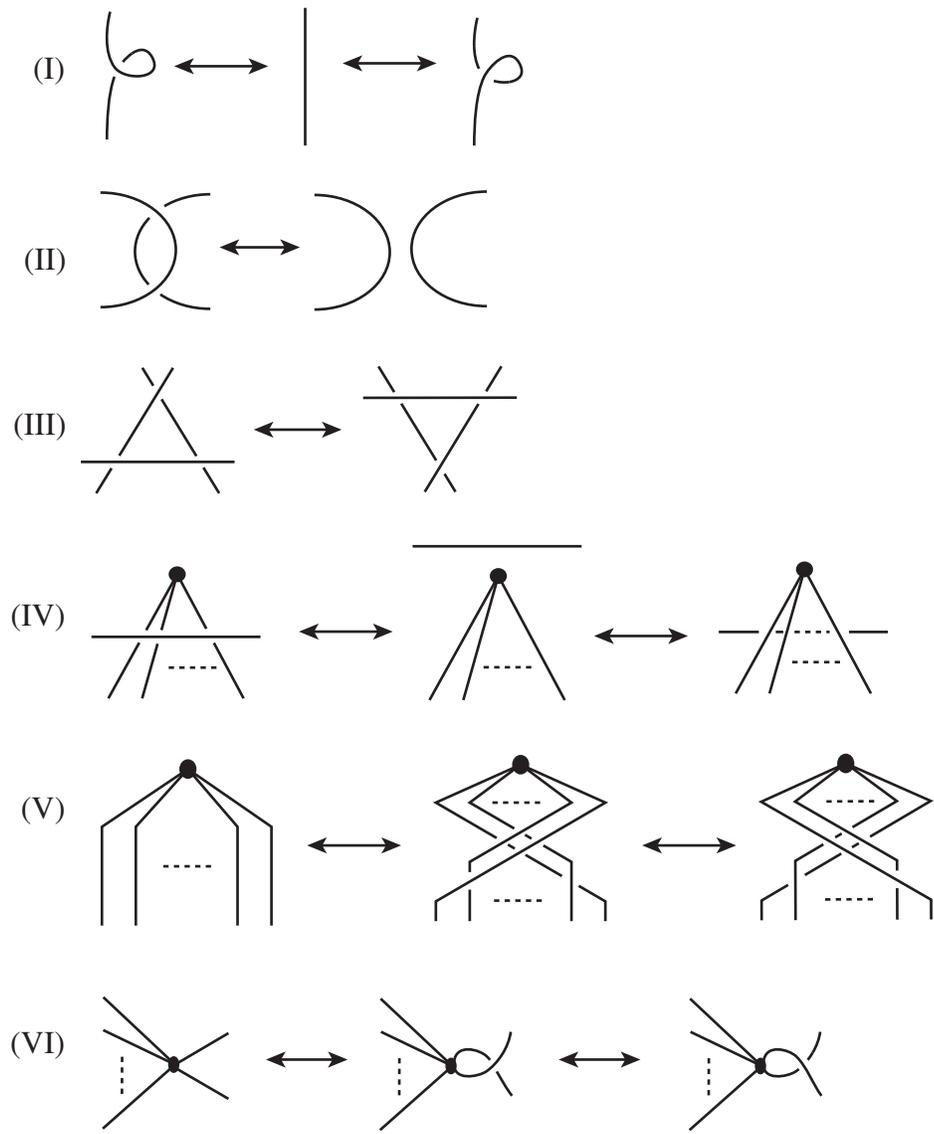}$$
    \caption{Reidemeister moves for graphs} \label{F:classicalmoves}
    \end{figure}

Our definition of virtual spatial graphs is combinatorial, and closely 
follows Kauffman's definition of virtual knots \cite{ka}.  First, we 
recall the definition of a classical spatial graph.  A {\it graph} is a 
pair $G = (V, E)$ of a set of {\it vertices V} and {\it edges} $E \subset 
V \times V$.  Unless otherwise stated, our graphs are connected and {\it 
directed}, so that each edge is an {\it ordered} pair of vertices.  A {\it spatial graph} is an 
embedding of $G$ in $\mathbb{R}^3$ that maps the vertices to points and an 
edge $(u, v)$ to an arc whose endpoints are the images of the vertices $u$ 
and $v$, and that is oriented from $u$ to $v$.  We will consider these 
embeddings modulo equivalence by ambient isotopy.  We can always represent 
such an embedding by projecting it to a plane so that each vertex 
neighborhood is a collection of rays with one end at the vertex and 
so that crossings of edges of the graph are transverse double points in the 
interior of the edges (as in the usual knot and link diagrams) \cite{ka2}.  
An example of such a diagram is shown in Figure~\ref{F:diagram}.
Kauffman \cite{ka2} and Yamada \cite{ya} have shown that ambient isotopy of spatial graphs is generated by a set of local moves on these diagrams which generalize the Reidemeister moves for knots and links.  These Reidemeister moves for graphs are shown in Figure~\ref{F:classicalmoves}.
The first five moves (moves (I) - (V)) generate {\it rigid vertex isotopy}, where the cyclic order of the edges around each vertex is fixed.  Moves (I) - (VI) generate {\it pliable vertex isotopy}, where 
the order of the vertices around each edge can be changed using move (VI).

\begin{figure} [ht]
    $$\includegraphics{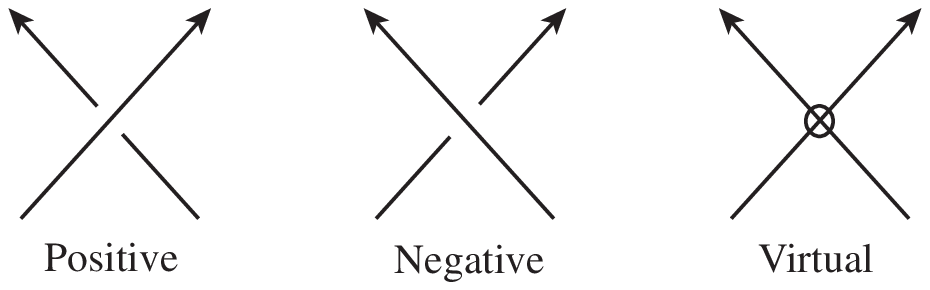}$$
    \caption{Types of crossings} \label{F:crossings}
    \end{figure}

\begin{figure} 
    $$\includegraphics{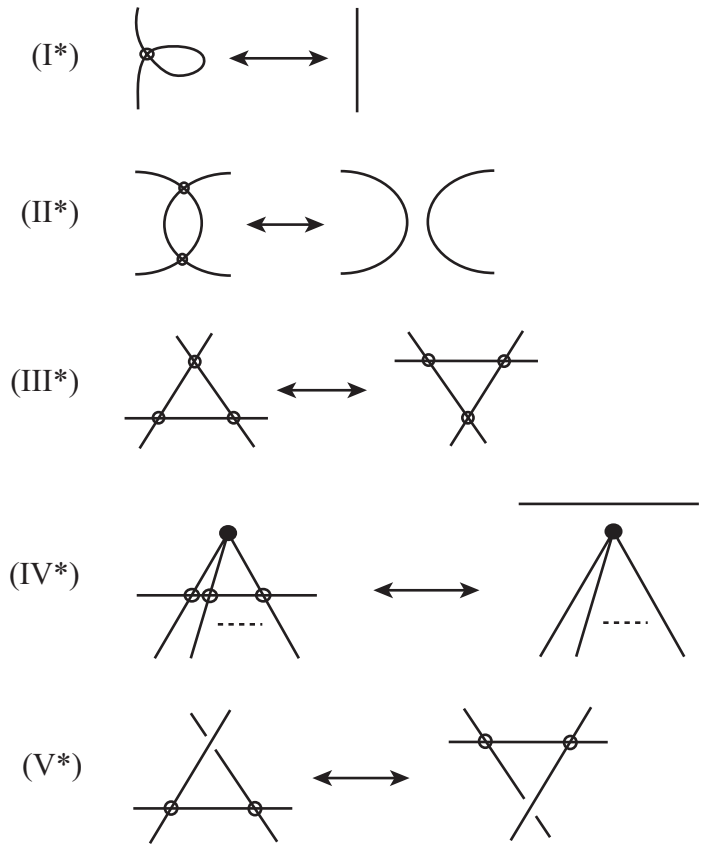}$$
    \caption{Reidemeister moves for virtual graphs} \label{F:virtualmoves}
    \end{figure}

A {\it virtual graph diagram} is just like a classical graph diagram, with 
the addition of {\it virtual crossings}.  We will represent a virtual 
crossing as an intersection of two edges surrounded by a circle, with no 
under/over information.  So we now have three kinds of crossings:  
positive and negative classical crossings and virtual crossings (see 
Figure~\ref{F:crossings}).
The idea is that the virtual crossings are not really 
there (hence the name ``virtual").  To make sense of this, we extend our 
set of Reidemeister moves for graphs to include moves with virtual 
crossings.  We need to introduce 5 more moves, (I*) - (V*), shown in 
Figure~\ref{F:virtualmoves}.  Notice that moves (I*) - (IV*) are just the 
purely virtual versions of moves (I) - (IV); move (V*) is the only move 
which combines classical and virtual crossings (in fact, there are two 
versions of the move, since the classical crossing may be either positive 
or negative).

\begin{figure} 
    $$\includegraphics{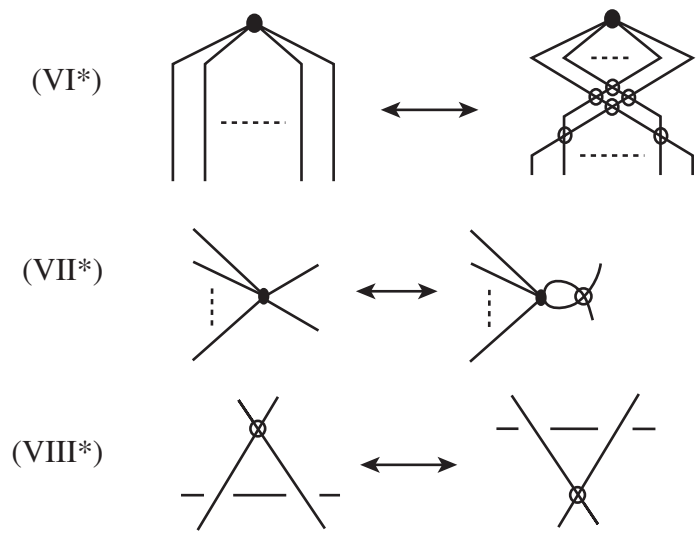}$$
    \caption{Forbidden Reidemeister moves for virtual graphs} \label{F:forbiddenmoves}
    \end{figure}

There are also three moves which, while they might seem reasonable, are {\it not} allowed.  These {\it forbidden moves} are shown in Figure~\ref{F:forbiddenmoves}.  We will explain why these moves are forbidden in the next section, when we discuss Gauss codes.

\subsection{Gauss Codes} \label{SS:gausscodes}

\begin{figure}
    $$\includegraphics{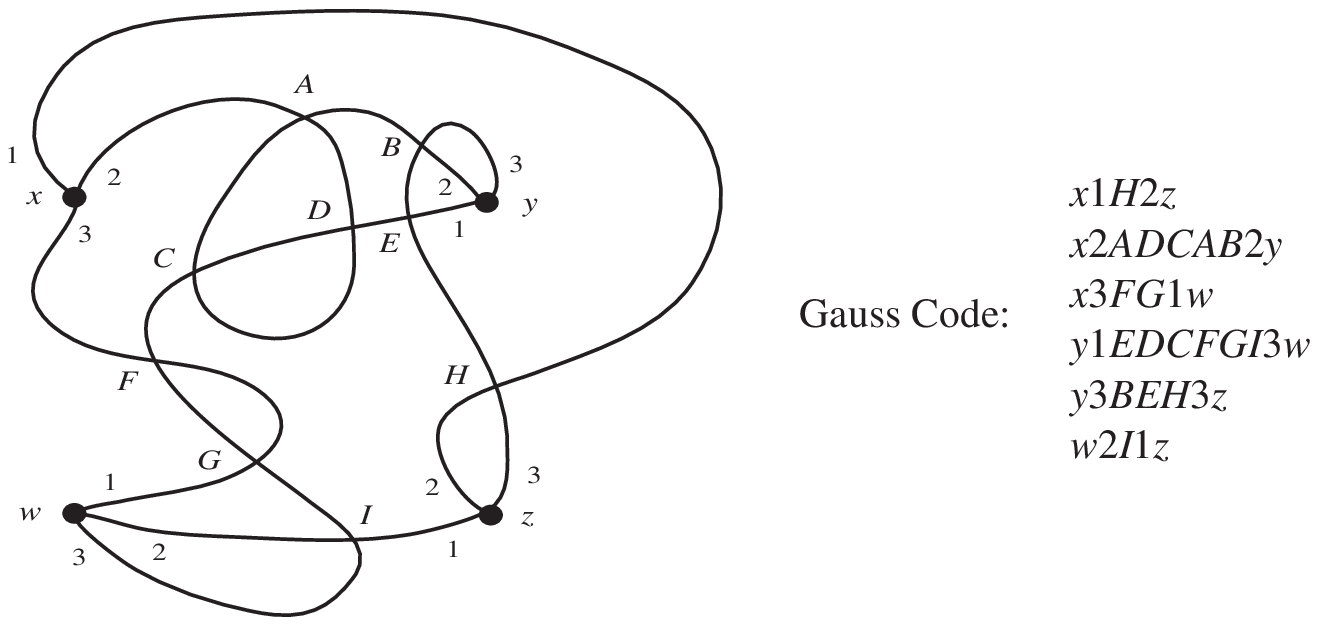}$$
    \caption{Gauss code for the shadow of a graph diagram} \label{F:gausscode}
    \end{figure}

\begin{figure}
    $$\includegraphics{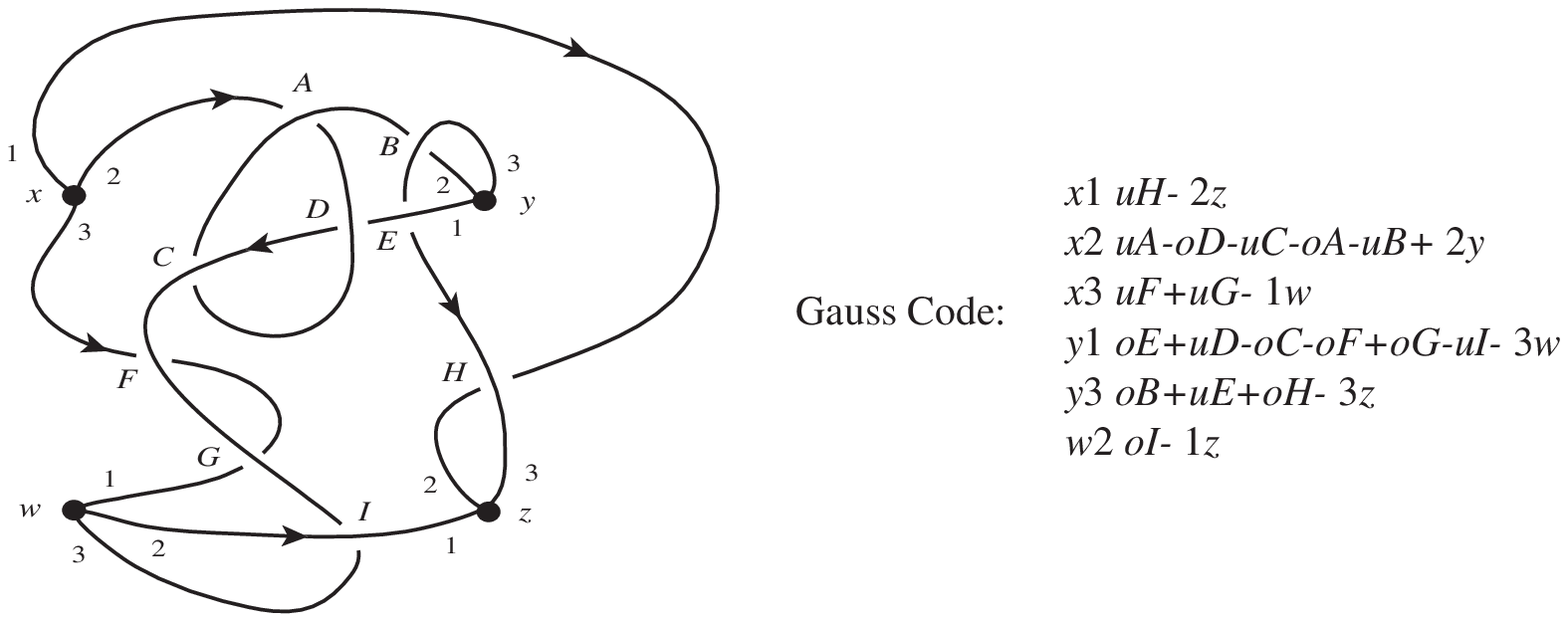}$$
    \caption{Gauss code for a graph diagram} \label{F:gausscode2}
    \end{figure}

One motivation for the study of virtual spatial graphs comes from the {\it Gauss code} of a diagram for a spatial graph.  Gauss codes for knot diagrams have a long history, and can be generalized to diagrams of any graph.  The Gauss code simply records the sequence of (labeled) crossings along each edge of the graph, so abstractly the Gauss code is just a set of sequences of symbols from some alphabet so that each symbol appears twice in the set.  Traditionally, the Gauss code is associated with an immersion of a closed curve (or graph) in the plane, so we begin by looking at the {\it shadow} of our graph diagram, where the over/under information at the crossings is ignored.  Figure~\ref{F:gausscode} illustrates how we write down the Gauss code for such a shadow.
To get the Gauss code for the original diagram, we can augment the Gauss code for its shadow by recording whether each crossing is an over-crossing (o) or an under-crossing (u).  If the graph is {\it directed} (i.e. the edges are oriented), we can also label each crossing by its sign, positive (+) or negative (-).  Figure~\ref{F:gausscode2} shows the Gauss code for a directed graph diagram.  (Clearly, a Gauss code must have the two occurrences of each symbol labeled with the same sign, and opposite over/under information.)

\begin{figure}
    $$\includegraphics{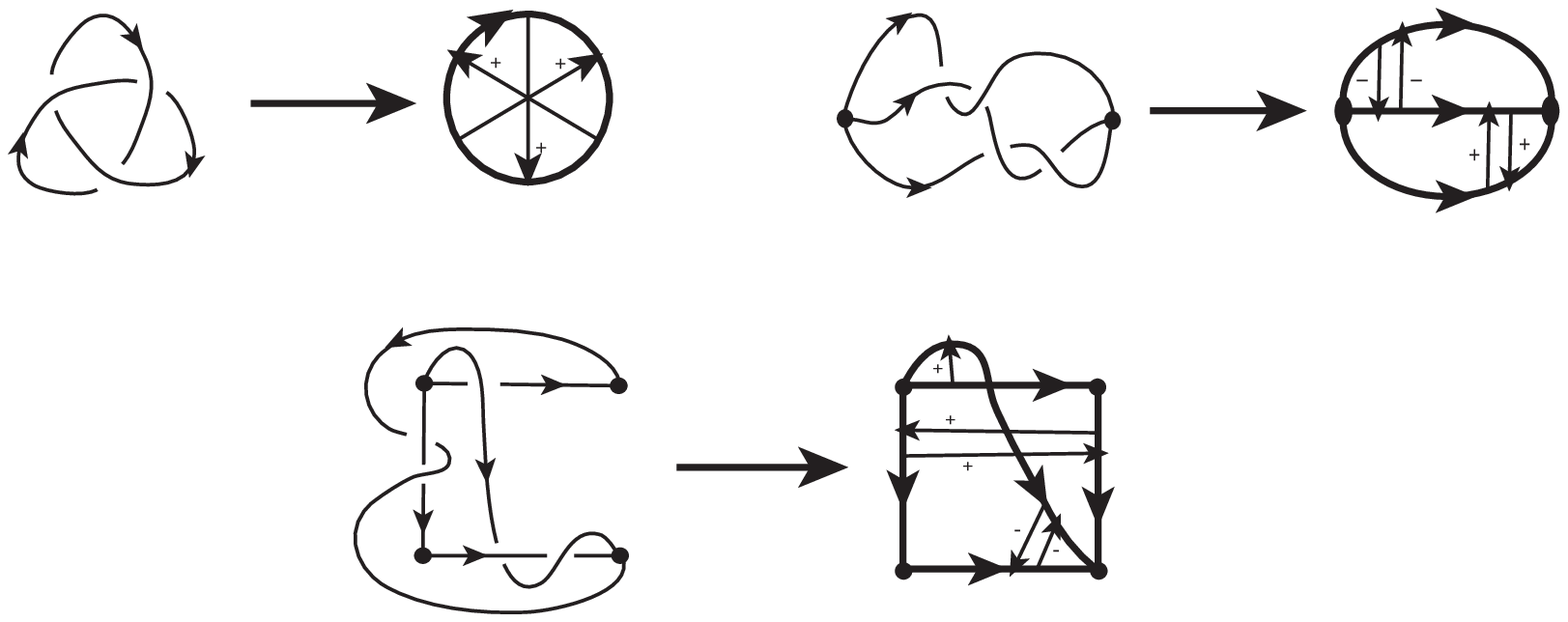}$$
    \caption{Examples of Gauss diagrams} \label{F:gaussdiagram}
    \end{figure}

It is often convenient to look at Gauss codes diagrammatically, using {\it Gauss diagrams} or {\it arrow diagrams}.  In these diagrams, we select a planar immersion of the graph, with the order of the edges at each vertex specified by the Gauss code, label points along the edges according to the crossing sequences of the Gauss code, and then draw an arrow between each pair of occurrences of a label.  The arrow is oriented from the undercrossing edge to the overcrossing edge, and is labeled with the sign of the crossing.  Once the arrows are drawn, the labels on the edges are redundant, and can be removed.  Some examples are shown in Figure~\ref{F:gaussdiagram}.

In the case of knots, there is an obvious immersion of the underlying graph - the standard embedding of a circle in the plane.  For spatial graphs, the situation is more complicated.  In general, there is no canonical choice for the ``simplest" immersion of the graph, though generally we will try to minimize the number of crossings.  But even if the underlying graph has a planar embedding, the order of the edges at the vertices may prevent us from using it, as in the third example in Figure \ref{F:gaussdiagram}.  So the Gauss diagram may have additional crossings, which we simply ignore.  Nevertheless, the Gauss diagram is a useful tool for understanding the effects of Reidemeister moves on Gauss codes.

\begin{figure}
    $$\includegraphics{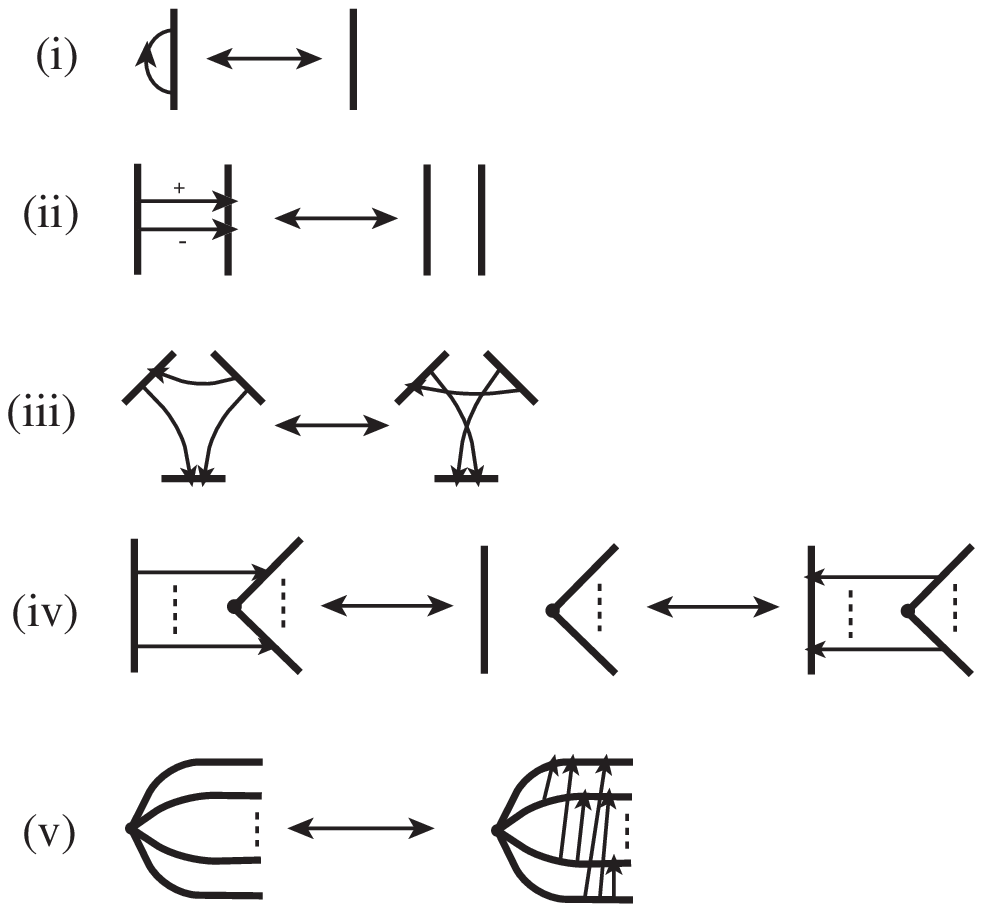}$$
    \caption{Reidemeister moves for Gauss diagrams} \label{F:gaussmoves}
    \end{figure}

The classical Reidemeister moves (I) - (V) for spatial graphs induce corresponding moves (i) - (v) on Gauss diagrams, shown in Figure~\ref{F:gaussmoves}.  So if spatial graphs $G_1$ and $G_2$ are equivalent modulo the Reidemeister moves (I) - (V), then their Gauss diagrams $D_1$ and $D_2$ are equivalent modulo moves (i) - (v).

However, the correspondence between equivalence classes of classical spatial graphs and equivalence classes of Gauss codes (or Gauss diagrams) is not bijective.  Even when we are considering only shadows, not all Gauss codes can be realized by the shadow of some graph diagram - those that can be are called {\it realizable} or {\it planar} codes.  An important problem in the study of Gauss codes is to find algorithms for determining whether a Gauss code is realizable.  For closed curves, there are several algorithms \cite{fo, de, rr}; the authors have generalized these methods to general graphs \cite{fm}.

Virtual graph diagrams also have Gauss codes, produced in exactly the same way, except that virtual crossings are ignored (hence we say that they are not ``real" crossings).  One motivation for studying virtual graph diagrams is that they allow us to realize the ``unrealizable" Gauss codes.

\begin{thm} \label{T:virtualrealize}
Every Gauss code can be realized as the code for a virtual graph diagram.
\end{thm}

{\sc Proof:}  We begin by embedding small neighborhoods of the vertices and the crossings in the plane, with the crossings decorated by orientation and over/under information.  The rest of the diagram consists of arcs between the vertices and crossings, and the Gauss code determines the endpoints and orientations of these arcs.  Simply draw in any collection of arcs with the desired endpoints, putting in virtual crossings wherever the arcs cross.  The result is a virtual graph diagram with the desired Gauss code.  $\Box$\\
\\
\noindent {\sc Remark:}  This proof will most likely not yield the ``best" virtual graph diagram.  It would be interesting (though undoubtedly difficult) to find an algorithm to produce a virtual graph diagram from a Gauss code with a minimal number of virtual crossings.  Presumably, finding this minimal number, the {\it virtual crossing number} would be as difficult as finding the classical crossing number of a knot or graph diagram.\\

So now we want to ask whether the correspondence between virtual graph diagrams (modulo moves (I) - (V) and (I*) - (V*)) and Gauss diagrams (modulo moves (i) - (v)) is a bijection.  In fact, we would like to {\it define} virtual spatial graphs as simply equivalence classes of abstract Gauss codes modulo moves (i) - (v).  However, while every virtual graph diagram has a well-defined Gauss code, it is possible for different diagrams to have the same Gauss code, so the inverse map may not be well-defined.  To show that this definition of a virtual spatial graph is the same as our original one, we need to show that two virtual graphs with the same Gauss code are virtually equivalent.\\
\\
\noindent {\sc Remark:} We can now see why the moves in Figure \ref{F:forbiddenmoves} are forbidden.  While moves (I*) - (V*) do not affect the Gauss code, moves (VI*) - (VIII*) {\it do} change the Gauss code.  Moves (VI*) and (VII*) change the order of the edges around a vertex, while move (VIII*) changes the order of two crossings along an edge.  Allowing these moves would force us to expand our list of moves on Gauss codes.

\begin{thm} \label{T:virtualgauss}
If two virtual graph diagrams have the same Gauss code, then they are virtually equivalent.
\end{thm}

{\sc Proof:}  Our proof follows the argument of the analogous theorem for virtual knots \cite{ka}.  Assume that $D$ and $E$ are two virtual graph diagrams with the same Gauss codes, so they have the same classical crossings, with the same local orientations and over/under behavior.  By an isotopy, we can assume that these classical crossings are in the same positions in the plane, and that a small neighborhood of the crossings is the same in both diagrams.  The arcs connecting these real crossings in each diagram contain only virtual crossings.  Say that the crossings $a$ and $b$ are connected by the arc $\gamma$ in $D$ and $\delta$ in $E$.  Since $\gamma$ and $\delta$ have the same endpoints, and have only virtual crossings, $\gamma$ may be moved to $\delta$ by virtual local moves (i.e. moves (I*) - (V*)) which do not change the Gauss code.  Doing this in turn with every arc in $D$ results in a virtual equivalence between $D$ and $E$, and completes the proof.  $\Box$

\begin{cor} \label{C:virtualclassical}
If a virtual graph diagram has a Gauss code which is realizable, then it is virtually equivalent to a classical graph diagram.
\end{cor}

We have not addressed the important question of whether the theory of classical spatial graphs is a proper subset of the theory of virtual spatial graphs, though we conjecture that it is.

\begin{conj}
If two classical spatial graphs are virtually equivalent, then they are classically equivalent.
\end{conj}

The corresponding result for virtual knots depends on the fact that the knot quandle (plus longitudes) is a complete knot invariant \cite{ka}.  While the quandle can be generalized to spatial graphs (both virtual and classical, as we will see in the next section), it is not nearly so powerful.  So it is not yet clear how to prove this conjecture in the broader context of spatial graphs.

\subsection{Forbidden Moves} \label{SS:forbidden}

If we allow the forbidden moves in Figure \ref{F:forbiddenmoves}, then many more virtual graph diagrams become equivalent.  In the case of knots, allowing move (VIII*) trivializes the theory, and all virtual knots become trivial \cite{ka, ne}.  However, when we look at virtual links or virtual graph diagrams, the effect is not quite so drastic - there are some properties which are not trivialized by the forbidden moves (see Proposition \ref{P:linkingnumber} in Section \ref{SS:T(G)}).  In this section we will briefly comment on the effects of allowing the forbidden moves.

\begin{prop} \label{P:forbidden}
Say that $G$ is a virtual graph diagram with Gauss diagram $D$.  Order the edges of $G$, and let $A_{i,j}$ be the set of arrows in $D$ with their tail on edge $i$ and their head on edge $j$.  If we allow move (VIII*) along with the moves for virtual rigid vertex isotopy (moves (I) - (V) and (I*) - (V*)), then we can transform $G$ (and hence $D$) such that:
\begin{enumerate}
	\item Every arrow of $D$ has its endpoints on distinct edges of the graph, so $A_{i,i}$ is empty.
	\item The arrows in $A_{i,j}$ are parallel and adjacent along both edges $i$ and $j$, and all have the same sign.
\end{enumerate}
If we also allow moves (VI), (VI*) and (VII*) (so we allow pliable vertex isotopy), then we can also ensure that: 
\begin{enumerate}
	\item $D$ has any desired ordering of the edges at each vertex.
	\item $A_{i,j}$ is empty whenever edges $i$ and $j$ are adjacent (i.e. are incident to the same vertex).
\end{enumerate}
\end{prop}
{\sc Proof:}  First, we will consider the effect of allowing the forbidden move (VIII*).  Nelson \cite{ne} proves that, in terms of Gauss diagrams, the move (VIII*) allows us move the head or tail of an arrow past the head or tail of an adjacent arrow (on the same edge).  This means that any arrow in the Gauss diagram of a graph which has its head and tail on the same edge can be moved until the head and tail are adjacent, and then erased by move (i).  So we can transform the Gauss diagram until $A_{i,i}$ is empty for every edge $i$.

Similarly, moving arrows past each other allows us to make all the arrows in $A_{i,j}$ parallel and adjacent.  If any pair of arrow in this set have opposite signs, then they cancel each other by move (ii).  So the remaining arrows all have the same sign. 

If, in addition, we allow move (VII*) (and therefore move (VI*)), then we can transpose two neighboring edges around any vertex.  If two edges $i$ and $j$ are adjacent at vertex $v$, then we can apply move (VII*) so that they are consecutive in the ordering of the edges around $v$.  Using move (VIII*), we can move any arrows in $A_{i,j}$ until there are no other arrows between their endpoints and vertex $v$.  Then we can use move (VI) (pliable vertex isotopy) to remove these arrows, leaving $A_{i,j}$ empty.  After removing all these edges, we can once again use move (VII*) to get any desired ordering of the edges around each vertex. $\Box$

\begin{quest} \label{Q:forbidden}
What is the equivalence relation on virtual graph diagrams (or their Gauss diagrams) generated by allowing move (VI*), but not move (VII*)?  What if we allow (VI*) and (VIII*)?
\end{quest}

\section{Invariants of Virtual Spatial Graphs} \label{S:invariants}

\subsection{The collection of virtual knots and links T(G)} \label{SS:T(G)}

In view of the large body of work on invariants of knots and links, a natural place to begin looking for invariants of spatial graphs, or virtual spatial graphs, is among the knots and links contained within the graph diagram.  In particular, we are interested in whether individual cycles of the graph are knotted, and whether {\it disjoint} cycles (cycles which do not share any edges or vertices) are linked.

\begin{figure} 
    $$\includegraphics{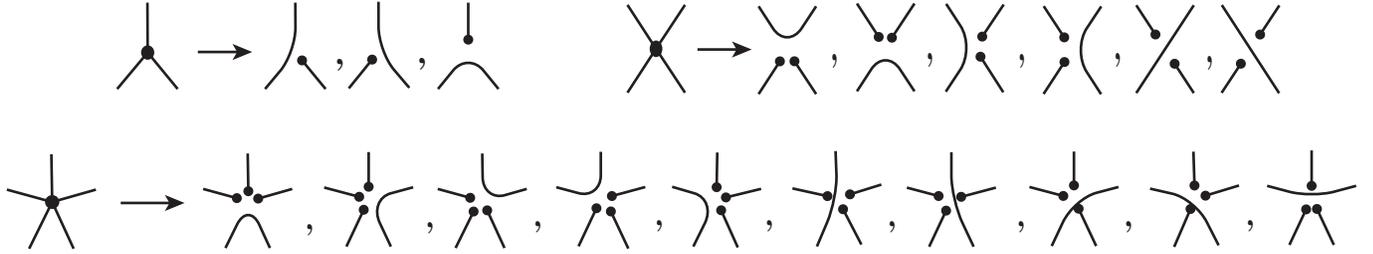}$$
    \caption{Local replacements of vertices} \label{F:replacement}
    \end{figure}

To formalize this, Kauffman \cite{ka2} introduced a topological invariant of a spatial graph (i.e. an invariant of pliable vertex isotopy) defined as the collection of all knots and links formed by a local replacement at each vertex of the graph.  Each local replacement joins two of the edges incident to the vertex and leaves the other edges as free ends (i.e. creates new vertices of degree one at the end of each of the other edges).  Figure~\ref{F:replacement} shows the possible replacements for vertices of degree 3, 4 and 5.
Choosing a replacement at each vertex of a graph $G$ creates a link $L(G)$ (after erasing all unknotted arcs).  $T(G)$ is the collection of all links $L(G)$ for all possible choices of replacements.  Kauffman showed that $T(G)$ is a pliable vertex isotopy invariant of $G$ \cite{ka2}.

\begin{figure} 
    $$\includegraphics{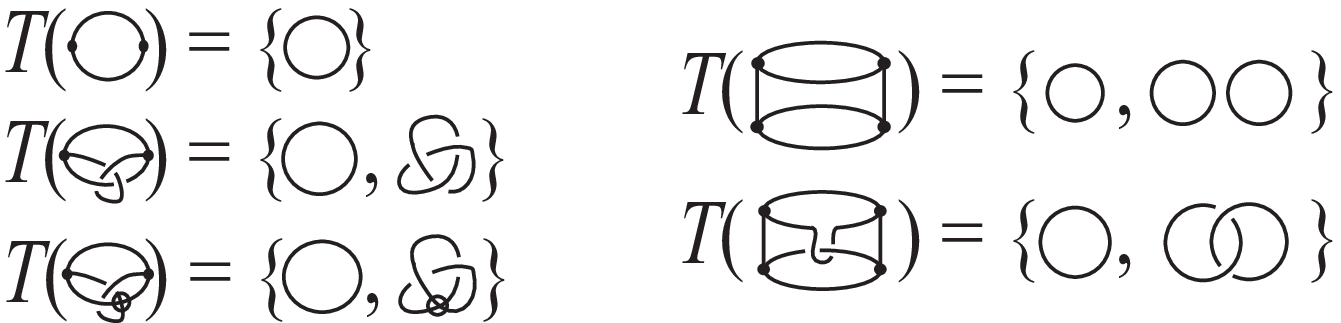}$$
    \caption{Examples of T(G) for virtual graph diagrams} \label{F:Texample}
    \end{figure}

For virtual graphs, we can define $T(G)$ in exactly the same way, except that it is now a collection of {\it virtual} links.  Kauffman's proof easily generalizes to show that $T(G)$ is also an invariant of virtual pliable vertex isotopy - in fact, $T(G)$ is invariant under all our moves except for move (VIII*) (including the forbidden moves (VI*) and (VII*)).  Figure~\ref{F:Texample} gives examples of $T(G)$ for some virtual spatial graphs.  It is worth observing that, since the trefoil knot with a single virtual crossing is not equivalent to a classical knot (as can be seen using the Jones polynomial), the third example on the left is a virtual spatial graph which is not equivalent to a classical spatial graph.

Once we have defined $T(G)$, we can obtain invariants of the virtual graph diagram by applying common knot and link invariants to the elements of $T(G)$.  For example, we can compute the pairwise linking numbers for all the links in $T(G)$ by using the Gauss formula ($\frac{1}{2}$(number of positive crossings) - $\frac{1}{2}$(number of negative crossings)).  If the links are virtual, these linking numbers may not be integers, but they are still invariant under all the classical and virtual Reidemeister moves.  In fact, the Gauss formula is also invariant under the forbidden moves (VI*), (VII*) and (VIII*).  For example, the two graphs on the right in Figure \ref{F:Texample} have links in $T(G)$ with different linking numbers, and so are inequivalent, even allowing the forbidden moves.  This shows that the forbidden moves do {\it not} trivialize virtual graph theory, as they do virtual knot theory.

\begin{prop} \label{P:linkingnumber}
There are virtual graph diagrams which are not equivalent modulo the forbidden moves.
\end{prop}

\subsection{Fundamental group} \label{SS:fundamental}

\begin{figure} 
    $$\includegraphics{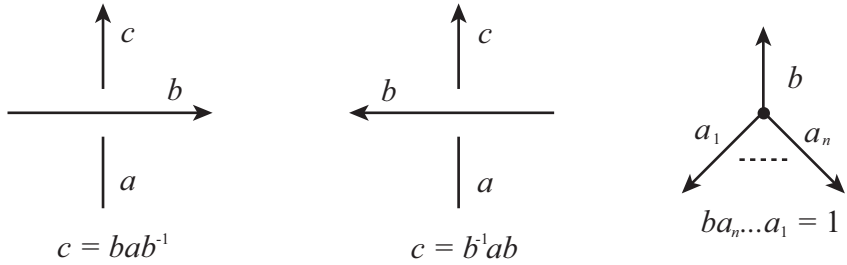}$$
    \caption{Wirtinger presentation of the fundamental group for a spatial graph} \label{F:wirtinger}
    \end{figure}

\begin{figure} 
    $$\includegraphics{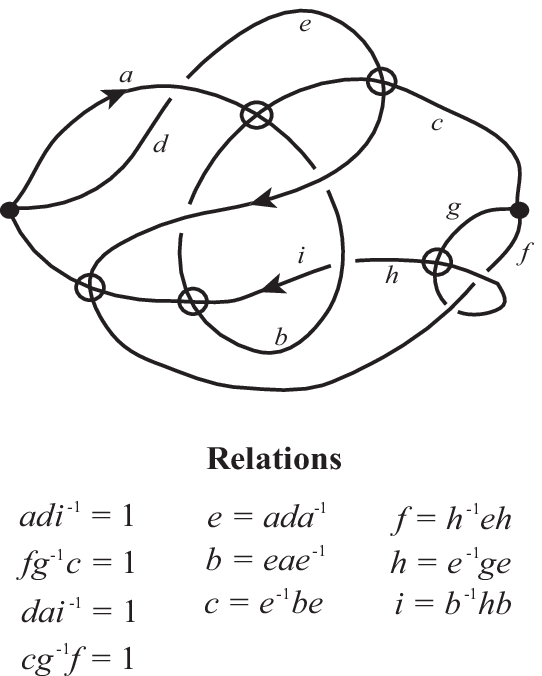}$$
    \caption{Fundamental group of a virtual spatial graph} \label{F:fgexample}
    \end{figure}

	The {\it fundamental group} of a classical knot or spatial graph is the fundamental group of its complement in $S^3$.  Given a diagram for the knot or graph, this group can be given a presentation, the {\it Wirtinger presentation}, involving one generator for each arc in the diagram and one relation for each crossing or vertex, as shown in Figure~\ref{F:wirtinger}.  In the relation at the crossing, changing the direction of an edge interchanges the corresponding generator in the word with its inverse.
Kauffman \cite{ka} defined the fundamental group of a virtual knot by constructing a presentation from a diagram using a generator for each arc between classical crossings, and a relation at each classical crossing.  We will define the fundamental group of a virtual spatial graph in the same way, by writing down a presentation with one generator for each arc between classical crossings (or vertices), and relations at each classical crossing or vertex (as shown in Figure \ref{F:wirtinger}).  An example is shown in Figure~\ref{F:fgexample}.
It is easy to check that the fundamental group is invariant under moves (I) - (VI) in Figure \ref{F:classicalmoves} (so it is an invariant of spatial graphs up to {\it pliable} vertex isotopy).  In fact, for a classical spatial graph, this is exactly the same as the classical fundamental group.  It is also easy to see that moves (I*) - (V*) in Figure \ref{F:virtualmoves} do not change any of the generators or relations, and so the fundamental group is an invariant of virtual spatial graphs.  However, all three of the forbidden moves in Figure \ref{F:forbiddenmoves} {\it do} change the fundamental group - in particular, moves (VI*) and (VII*) change the relation at a vertex.

\subsection{Virtual Graph Quandle} \label{quandle}

The quandle is a combinatorial knot invariant that was generalized to 
virtual knots by Kauffman \cite{ka}, and strengthed by Manturov \cite{ma}. 
Modifying Manturov's approach, we can construct a similar invariant for 
virtual spatial graphs, though for general graphs this invariant is less 
potent than in the case of knots.

Let $M$ be a set with one symbol for each arc in a diagram of $G$. 
Further, let $M$ have an operation $\circ$, an involution $a \rightarrow 
\overline{a}$, and an invertible function $f$.  Construct the set $X$ of 
all words in the elements of $M$ using $\circ$, $\overline{a}$ and $f$.

\begin{figure} 
    $$\includegraphics{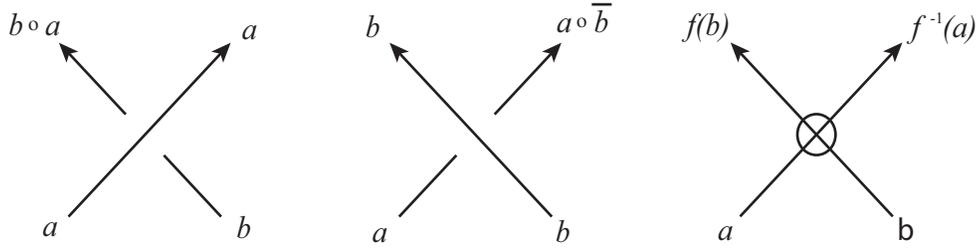}$$
    \caption{The virtual quandle crossing relations.} \label{crossingrel}
    \end{figure}

The virtual graph quandle $Q(G)$, an invariant of the virtual spatial 
graph $G$, is formed from $X$ by quotienting out the relations listed 
below.  Following each relation, we note the Reidemeister move(s) that 
require that relation.  To encode information about the diagram, we 
identify edge labels so that the arcs meeting in a crossing are labeled as 
in Figure \ref{crossingrel}.

\begin{figure} 
    $$\includegraphics{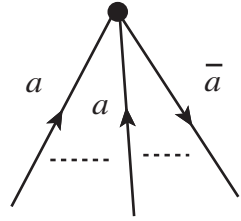}$$
    \caption{The virtual quandle relation at a vertex.} \label{virtvert}
    \end{figure}

In addition, we require that all arcs meeting in a vertex be labeled as in 
Figure \ref{virtvert}. That is, we identify the labels of two arcs 
entering the vertex, and we identify the label of an arc entering the 
vertex with the bar of the label of an arc leaving the vertex.  This 
restriction at the vertices is needed to ensure invariance under move 
(IV).

The first three relations are necessary for invariance under Reidemeister moves on the edges of $G$, and are the same as the relations needed for the knot quandle.

\begin{align*}
&a \circ a = a  &\mbox{(I)}\\
&(a \circ b) \circ \overline{b} = (a \circ \overline{b}) \circ b = a  
&\mbox{(II)}\\
&(a \circ b) \circ c = (a \circ c) \circ (b \circ c)  &\mbox{(III)}\\
\end{align*}

We now need to add relations to give invariance under classical Reidemeister moves involving a vertex.  Let $d =$ gcd($d_{1}, d_{2}, \ldots d_{k}$), where $d_{i}$ is the valence of vertex $i$.

\begin{align*}
&((b \circ a) \circ a) \ldots \circ a) = ((b \circ 
\overline{a}) \ldots \circ \overline{a}) = b \quad \mbox{where $a$ 
occurs $d$ times} &\mbox{(IV)}\\
&\overline{(a \circ b)} = \overline{a} \circ b &\mbox{(IV)}\\
&\overline{a} \circ a = \overline{a} &\mbox{(V),(VI)}\\
&\overline{(\overline{a})} = a &\mbox{(V),(VI)}\\
\end{align*}

The following relations ensure invariance under the virtual Reidemeister 
moves.  Notice that due to the labeling rule of Figure \ref{crossingrel}, 
$Q(G)$ is already invariant under moves (I*), (II*), and (III*).

\begin{align*}
&f^{d}(b) = b &\mbox{(IV*)}\\
&\overline{f(a)} = f(\overline{a}), \overline{f^{-1}(a)} = 
f^{-1}(\overline{a}) &\mbox{(V*)}\\
&f(a \circ b) = f(a) \circ f(b) &\mbox{(V*)}\\
\end{align*}

\begin{figure}
    $$\includegraphics{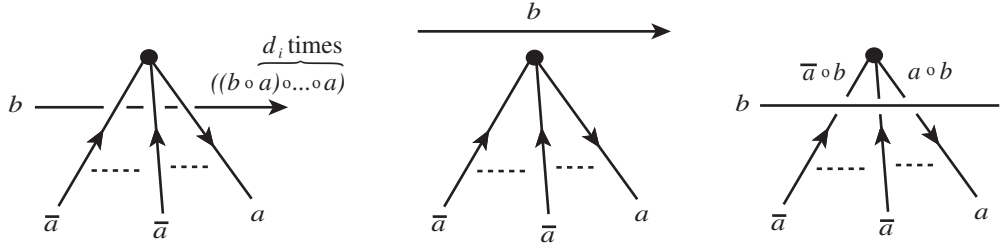}$$
    \caption{Invariance under IV. The argument is similar for other choices of orientation.} \label{r4check}
    \end{figure}

\begin{figure}
    $$\includegraphics{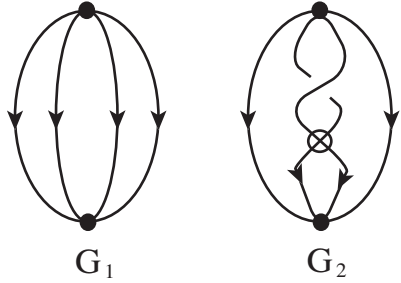}$$
    \caption{$Q(G_{1})$ is the free quandle on $a$ with $d=4$.  $Q(G_{2})$ is the free quandle on $a$ with $d=4$ and $f(a) = a$.} \label{quaneg}
    \end{figure}

Figure \ref{r4check} demonstrates the invariance of $Q(G)$ under Reidemeister move IV. For the left side, we use $((b \circ a) \circ a) \ldots \circ a) = b$ to show invariance, and for the right, $\overline{a} \circ b = \overline{(a \circ b)}$ to show that the vertex relation is still satisfied.  Invariance under the other Reidemeister moves can be checked in a similar manner.  Figure \ref{quaneg} depicts two virtual spatial graphs that are distinguished by their quandles.  

\subsection{Yamada polynomial} \label{SS:yamada}

Yamada introduced a polynomial invariant $R$ of spatial graphs in \cite{ya}.  In this section we will review the definition of this invariant, and show that it can be extended to an invariant of virtual spatial graphs.

Yamada's polynomial for an undirected graph $G$ can be defined combinatorially using skein relations as the unique polynomial $R(G)(A)$ which satisfies the following formulas:

\begin{enumerate}
    \item $R(\includegraphics{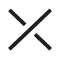}) = AR(\includegraphics{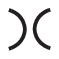})+A^{-1}R(\includegraphics{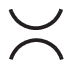}) + R(\includegraphics{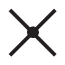})$
    \item $R(G) = R(G\backslash e) + R(G/e)$, where $e$ is a nonloop edge in $G$, $G \backslash e$ is the result of deleting $e$, and $G/e$ is the result of contracting $e$.
    \item $R(G_1 \amalg G_2) = R(G_1)R(G_2)$, where $\amalg$ denotes disjoint union.
    \item $R(G_1 \vee G_2) = -R(G_1)R(G_2)$, where $G_1 \vee G_2$ is the graph obtained by joining $G_1$ and $G_2$ at any single vertex.
    \item $R(B_n) = -(-\sigma)^n$, where $B_n$ is the $n$-leafed bouquet of circles and $\sigma = A + 1 + A^{-1}$.  In particular, if $G$ is a single vertex, $R(G) = R(B_0) = -1$.
    \item $R(\emptyset) = 1$
\end{enumerate}

Using these skein relations, $R(G)$ can be computed by reducing the graph $G$ to a bouquet of circles.  $R(G)$ is an invariant of spatial graphs up to {\it regular rigid vertex isotopy}, meaning that it is invariant under moves (II), (III) and (IV) in Figure \ref{F:classicalmoves}, but not moves (I), (V) or (VI) \cite{ya}.  The behavior of $R$ under these moves are shown in Figure~\ref{F:yamada} (see \cite{ya}):

    \begin{figure}
    $$\includegraphics{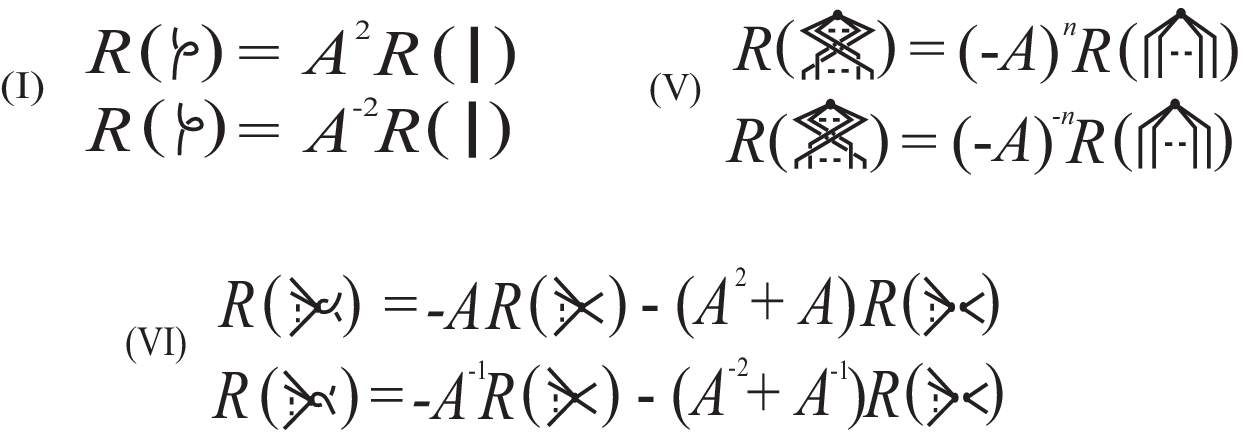}$$
    \caption{Behavior of $R(G)$ under moves (I), (V) and (VI)} \label{F:yamada}
    \end{figure}

From these formulas we see that we can obtain an invariant of rigid vertex isotopy (invariance under moves (I) - (V)) by defining $\bar{R}(G) = (-A)^{-m}R(G)$, where $m$ is the smallest power of $A$ in $R(G)$.  This will still not be invariant under move (VI), however, so it is not an invariant of pliable vertex isotopy.

\begin{figure}
    $$\includegraphics{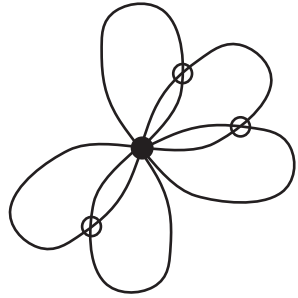}$$
    \caption{A virtual bouquet of circles} \label{F:virtualbouquet}
    \end{figure}

\begin{figure}
    $$\includegraphics{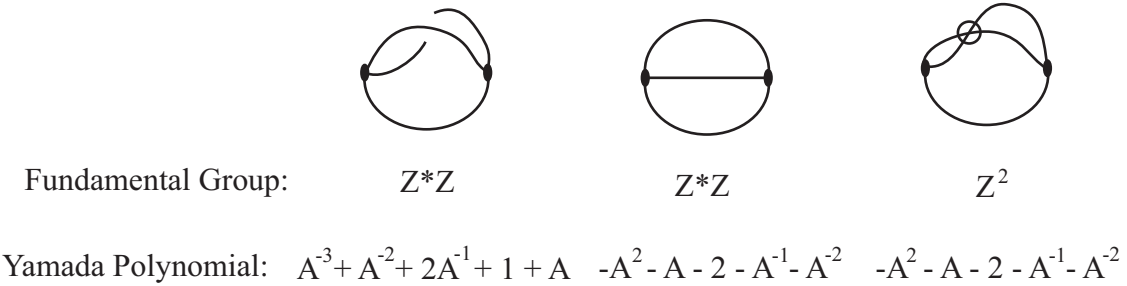}$$
    \caption{Diagrams for the $\theta$-graph} \label{F:examples}
    \end{figure}

In the case of virtual spatial graphs we can use exactly the same skein relations to compute $R(G)$ and $\bar{R}(G)$, simply by ignoring virtual crossings.  The only difference is that we may end up with a {\it virtual bouquet} - a bouquet of circles with only virtual crossings, as in Figure~\ref{F:virtualbouquet}.
We will simply ignore the virtual crossings.  In other words, if $G$ is a virtual bouquet of $n$ circles, then $R(G) = R(B_n) = -(-\sigma)^n$.  Given a planar graph, a planar diagram for the graph is called {\it trivial}.  It is an open question whether the Yamada polynomial for a non-trivial classical diagram of a planar graph can be the same as the polynomial for the trivial diagram.  However, it is easy to find nontrivial virtual graph diagrams of planar graphs with trivial Yamada polynomial.  Figure~\ref{F:examples} gives three diagrams of the $\theta$-graph, including the trivial one, together with their fundamental groups and Yamada polynomials.  These show that there are non-trivial diagrams with the same fundamental group as the trivial diagram but a different Yamada polynomial, and vice-versa.

\small

\normalsize

\end{document}